\documentclass{article}
\usepackage{latexsym}
\usepackage{epsfig}
\usepackage{amsmath}
\usepackage{amssymb}
\usepackage{amsthm}
\newtheorem{theorem}{Theorem}[section]

\theoremstyle{remark}

\theoremstyle{definition}
\newtheorem{definition}[theorem]{Definition}

\newcommand{\R}{{\mathbb R}} 

\begin{document}

\title{Optimal step length for the maximal decrease of a self-concordant function by the Newton method}

\author{Anastasia Ivanova\thanks{%
Univ.\ Grenoble Alpes, LJK, 38000 Grenoble, France; HSE University, Moscow, Russian Federation
({\tt anastasia.ivanova@univ-grenoble-alpes.fr}).}
\and
Roland Hildebrand \thanks{%
Univ.\ Grenoble Alpes, CNRS, Grenoble INP, LJK, 38000 Grenoble, France; MIPT, Institutskiy Pereulok 9, 141701 Dolgoprudny, Russian Federation
({\tt roland.hildebrand@univ-grenoble-alpes.fr}).}}

\maketitle

\begin{abstract}
    In this paper we consider the problem of finding the optimal step length for the Newton method on the class of self-concordant functions, with the decrease in function value as criterion. We formulate this problem as an optimal control problem and use optimal control theory to solve it.

\end{abstract}

\section{Introduction}

This paper is devoted to the problem of finding the optimal step length of Newton's method on the class of self-concordant functions, motivated by the appearance of this class in barrier methods for conic programming. Self-concordant functions were introduced by Yu. Nesterov and A. Nemirovsky \cite{nesterov1994interior} when studying the behavior of Newton's method, as follows.
\begin{definition}
A convex $C^3$ function $f \, : \, D \rightarrow \R$  on a convex domain $D$ is called self-concordant if it satisfies the inequality
\begin{equation}
\label{self_conc}
    |f'''(x)[h,h,h]| \leq 2(f''(x)[h,h])^{3/2}
\end{equation} 
for all $x \in D$ and all tangent vectors $h$.

It is called \textit{strongly self-concordant} if in addition $\lim_{x \rightarrow \partial D} f(x) = +\infty$. 
\end{definition}

Step lengths for the damped Newton method were also considered in  \cite{burdakov1980some, ralph1994global, nesterov2018lectures}.

In this paper we find the optimal step length of Newton's method with respect to the decrease of the function value. This criterion was considered in \cite[Theorem 2.2.1]{nesterov1994interior}, where the decrease has been lower bounded by an explicit function of the step length $\gamma_k$ and the Newton decrement $\rho_k$. The same bound has been derived in \cite{gao2019quasi} in a more general context. In the latter paper it is shown that the step length $\gamma_k = \frac{1}{1 + \rho_k}$ maximizes this lower bound. The same expression for the step length is also proposed in \cite[Theorem 2.2.3]{nesterov1994interior} for larger values of the decrement. In the present paper we show by employing optimal control theory that this step length is optimal. 

Optimal control theory has already been used in \cite{Hildebrand21Newton} to find an optimal step-length $\gamma^*$ which minimizes the worst-case Newton decrement in the next iteration.

The remainder of the paper is structured as follows. In Section \ref{Section:Problem_statement} we describe the problem statement.   In Section \ref{Section:Solve-problem} we rephrase the problem as an optimal control problem and solve it analytically.

\section{Problem statement}
\label{Section:Problem_statement}

We consider Newton's method with a damped step
\begin{equation}
\label{newton_step}
    x_{k+1} = x_k - \gamma_k (F''(x_k))^{-1} F'(x_k),
\end{equation}
where $\gamma_k \, \in \, (0,1]$ is the step-size.
Since Newton's method is affinely invariant, i.e. a sequence of iterations on a given function transforms into itself under an affine transformation of coordinates, it is natural to study the behavior of the method on a class of functions that is also affinely invariant. This leads to the self-concordant functions which naturally arise as an affinely invariant analogue of functions with a Lipschitz continuous Hessian, and hence is well suited for an analysis of the behaviour of Newton's method.

Namely \cite{nesterov1994interior}, the state at iteration $k$ is described by a single scalar, the \emph{Newton decrement}
\begin{equation}
\label{decrement}
    \rho_k = ||F'(x_k)||_{F''(x_k)} = \sqrt{(F'(x_k))^{\top}(F''(x_k))^{-1}F'(x_k)}.
\end{equation}

We consider the problem of finding a step length $\gamma_k$ which maximizes the decrease $F(x_{k+1})-F(x_k)$ of the function value in the worst case realization of the function $F(\cdot)$. So, we firstly need for given step length and given decrement to find the worst realization of the function giving the minimal decrease, and then to maximize this progress over the value of the step length, yielding the optimal step length as a function of the decrement. This leads to the following optimization problem:
\begin{equation} \label{minmax_problem}
   \max_{\gamma_k}\min_{F \in {\cal S}} \left(F(x_{k}) - F(x_{k+1})\right),
\end{equation}
where $\gamma_k$ is the step length, $x_{k+1}$ is given by \eqref{newton_step}, the decrement $||F'(x_k)||_{F''(x_k)}$ is fixed to some value $\rho_k$, and ${\cal S}$ is the class of functions satisfying \eqref{self_conc}.

\section{Solution using optimal control theory} 
\label{Section:Solve-problem}

In this section we describe the solution of problem \eqref{minmax_problem}. 

We consider a single iteration of the Newton method. Let the end point $x_{k+1}$ be given by \eqref{newton_step} and consider the line segment between $x_k$ and $x_{k+1}$. We study the evolution of the values of the function and its derivatives along this segment. This distinguishes our approach from the approach in \cite{de2020worst}, where $n$ iterations and only the values of the function and its derivatives at the points $x_{1}, \ldots, x_{n}$, i.e., a finite dimensional object, are considered. In contrast to this we consider an infinite dimensional object.  The suitable apparatus to solve this problem is optimal control theory. We start with formulating our problem as an optimal control one. 

\subsection{Optimal control problem statement}

For simplicity we denote $\gamma := \gamma_k$  and $\rho: = \rho_k$. %First of all, we may make some simplifications due to the affine invariance of the problem setup. So, we move the current iterate to the origin, $x_k = 0$, normalize the Hessian at the initial point to $ F''(x_k) = I$ by scaling the coordinate $x$, and possibly flip the real axis to achieve $F'(x_k) = -\rho e_1$. Then we get $ x_{k+1} = \gamma \rho e_1.$
Parameterize the line segment between the current iterate $x_k$ and the next one $x_{k+1}$ affinely by a variable $s \in [0,T]$, such that 
\begin{equation}
\label{param_step}
    x(0) = x_k, \quad x(s) = x_k + sh, \quad x(T) = x_{k+1} = x_k + Th,
\end{equation}
where $h$ is a direction vector along the segment between $x_k$ and $x_{k+1}$. Later we shall impose a normalization condition on the vector $h$, thereby determining the value of $T$.

Let us investigate the evolution of the gradient and Hessian of $F$ along the segment. %For simplicity we introduce the functions $g(s) = F'(\sigma(s))$ and $H(s) = F''(\sigma(s))$.
Since our goal is to find a realization of the function $F(\cdot)$, which minimizes the decrease of the function value, we represent the cost function  as follows: 
\begin{eqnarray*}
  &&  F(x_{k}) - F(x_{k+1}) = - \int_{0}^{T} \langle F'(x_k + s h), h \rangle\, ds.
  %= \int_{s_0}^{s_1} \langle g(s), \dot x_1(s) e_1 \rangle ds = \int_{s_0}^{s_1}  \dot x_1(s)  g_1(s)  ds.
\end{eqnarray*}
We get the maximization problem 
\begin{eqnarray*}
  &&  \int_{0}^{T} \langle F'(x_k + s h), h \rangle \,ds \rightarrow \max.
\end{eqnarray*}
For simplicity we introduce the function $g(s) =  \langle F'(x_k + s h), h \rangle $, then 
\begin{equation*}
    \xi(s) := \frac{dg(s)}{ds}= \langle F''(x_k + s h)h, h \rangle.
\end{equation*}
From \eqref{self_conc}  we get 
\begin{equation*}
    \left |\frac{d \xi(s)}{d  s} \right | = |F'''(x_k + sh)[h, h,h]| \leq 2 ( F''(x_k + s h)[h, h ] )^{3/2} = 2\xi(s)^{3/2},
\end{equation*}
or equivalently 
\begin{equation}
\label{eq_contr_sys}
   \frac{d \xi(s)}{d  s}  = 2u\xi(s)^{3/2},
\end{equation}
where $u \in [-1, 1]$. Thus, $u$ can be interpreted as a control, $U = [-1, 1]$ as the set of admissible controls, and \eqref{eq_contr_sys} as a controlled system, where $\xi(s)$ is a positive scalar function.
Introducing the function  $w(s) := \sqrt{\xi(s)}$, we get from \eqref{eq_contr_sys}  that
\begin{equation*}
   \frac{d w^2(s)}{d  s}  = 2u w (s)^{3} \quad \Rightarrow \quad  \frac{d w(s)}{d  s}  = u w (s)^{2}.
\end{equation*}
To impose a normalization condition on $h$, note that 
\begin{equation*}
    ||h||_{F''(x_k)} := \sqrt{F''(x_k)[h, h]} = w(0),
\end{equation*}
since $x_k$ corresponds to the value $s=0$. Moreover, from \eqref{newton_step} and \eqref{param_step} we get 
\begin{equation}
\label{eq_for_T_1}
    x_{k+1} - x_k = Th = -\gamma (F''(x_k))^{-1}F'(x_k).
\end{equation}
Then multiplying by $F'(x_k)$ and using \eqref{decrement}, we get 
\begin{equation}
\label{eq_for_T_2}
    T \langle F'(x_k), h \rangle = - \gamma \rho^2.
\end{equation}
Moreover, multiplying \eqref{eq_for_T_1} by $F''(x_k)$ and $h$ , we obtain
\begin{equation}
\label{eq_for_T_3}
    T F''(x_k)[h, h] = - \gamma \langle F'(x_k), h \rangle.
\end{equation}
Substituting this in \eqref{eq_for_T_2}, we get 
\begin{equation*}
    T^2 F''(x_k)[h, h] = \gamma^2 \rho^2.
\end{equation*}
Choosing the normalization for $h$ such that $||h||_{F''(x_k)} = 1$, we obtain that $T = \gamma \rho$. Substituting this into \eqref{eq_for_T_3}, we get that at $s=0$
\begin{equation*}
    w(0) = 1, \quad g(0) = - \rho.
\end{equation*}

Finally, we get the following optimal control problem 
\begin{eqnarray*}
&&  \int_{0}^{\gamma \rho } g(s) d s \, \rightarrow \, \text{max},\\
&& \frac{dg}{ds} = w^2, \quad \frac{dw}{ds} = uw^2,\\
&& w(0) = 1, \quad g(0) = -\rho
\end{eqnarray*}
with control $u \in [-1,1]$. Introducing a new variable $t$, such that 
\begin{eqnarray*}
t = - ||x_{k+1} - x||_{F''(x)} = - \sqrt{\xi(s)[h, h] \cdot ( T - s)^2} = -w \cdot (T - s),
\end{eqnarray*}
we get 
\begin{eqnarray*}
\frac{dt}{ds} &= &-w^2 u \cdot (T - s) + w = w \cdot (ut + 1),\\
\frac{dg}{d t} &= &\frac{dg}{d s}\cdot \frac{ds}{dt}= \frac{w}{1+ut},\\
\frac{dw}{d t} &= &\frac{dw}{d s}\cdot \frac{ds}{dt}= \frac{wu}{1+ut},
\end{eqnarray*}
and $t \in [-\rho \gamma, 0]$. Denoting $z = \frac{g}{w}$, we obtain 
\begin{eqnarray*}
\frac{d z}{d t } &=& \frac{d (g/w)}{d t } = \frac{\tfrac{w}{1 + ut}\cdot w - \tfrac{g}{w}\cdot \tfrac{w^2 u}{1 + ut}}{w^2} = \frac{1 - z u}{1 + ut}, \\
g \,ds &= &\frac{g}{w(1 + ut)}dt = \frac{z}{1 + ut} dt.
\end{eqnarray*}
So, we can rewrite the  optimal control problem as follows 
\begin{eqnarray}
\label{final_problem}
&&  \int^0_{-\gamma \rho} \frac{z}{1+ut} dt \, \rightarrow \, \text{max},\\
&& \dot z = \frac{1-uz}{1+ut}, \notag\\
 &&  z(-\gamma \rho) = -\rho, \notag
\end{eqnarray}
where $u \, \in \, U = [-1, 1]$ and $z \, \in \, \R$.

\subsection{Solution of the problem}

In this section we solve problem \eqref{final_problem}. 

Firstly, according to Pontryagin's maximum principle  \cite{pontr}, we get the following Hamiltonian for the optimal control problem \eqref{final_problem} 
\begin{equation*}
    \mathcal{H} = \frac{z}{1+ut} + \psi \frac{1-uz}{1+ut},
\end{equation*}
where $\psi \in \mathbb{R}$ is the adjoint variable to $z$. The dynamics of the adjoint variable is given by
\begin{eqnarray*}
 \dot \psi = - \frac{\partial \mathcal{H}}{\partial z}  =  \frac{u\psi - 1}{1+ut} .
\end{eqnarray*}

If the control is bang-bang, i.e., $u$ is piece-wise constant with values in $\{-1, 1\}$, then the dynamics of the primal and adjoint variables can be integrated explicitly. For the primal variable $z$  we get 
\begin{equation*}
    -u\log|1-uz| + C' = u\log|1+ut| \quad \Rightarrow \quad z(t) = \frac{C + t}{1 + ut},
\end{equation*}
where $C = (1 \pm e^{C'u})u$ is some constant.  For the adjoint variable $\psi$ we get the solutions
\begin{eqnarray*}
u\log|u\psi-1| + C'' = u\log|1+ut| \quad \Rightarrow \quad \psi(t) = C_1(tu+1) +  \frac{1}{u},
\end{eqnarray*}
where $C_1 = \pm e^{-C''u}$ is some constant. The transversality condition at the end-point gives $\psi(0) = 0$, hence $C_1 = -\frac{1}{u}$ and 
\begin{equation*}
    \psi(t) = -t
\end{equation*}
for all $t$ from the last switching point to $t = 0$. 

Our next step is to determine the optimal control $u$. Maximizing the Hamiltonian over the variable $u$ we get 
\begin{equation*}
u = 
  \begin{cases} 
   1, & \text{if } \psi z + tz + t\psi < 0, \\
   -1,       & \text{if } \psi z + tz + t\psi > 0.
  \end{cases}
\end{equation*}
For $t$ sufficiently close to 0 we get
\begin{equation*}
    \psi z + tz + t\psi  = -t\cdot  \frac{C + t}{1 + ut} + \left(  \frac{C + t}{1 + ut} - t \right)\cdot t = -t^2 < 0 ,
\end{equation*}
and the control  $u = 1$ is optimal. But the expression $-t^2$ remains negative for all negative $t$ up to the starting point $t = -\gamma\rho$. Hence the control $u = 1$ and above expressions for $z(t),\psi(t)$ are valid over the whole interval.

Using the boundary conditions for  $z(t)$, we obtain that 
\begin{equation*}
    z(-\gamma \rho ) = \frac{C - \gamma \rho}{1  - \gamma \rho} = -\rho \quad \Rightarrow \quad C = -\rho +\gamma\rho^2 + \gamma\rho.
\end{equation*}
Therefore
\begin{equation*}
    z(t) = \frac{-\rho +\gamma\rho^2 + \gamma\rho + t}{1  + t}.
\end{equation*}
Substituting this value into the objective function, we obtain 
\begin{eqnarray*}
    - \int^{0}_{-\gamma \rho} \frac{z(t)}{1+t} \,dt &= &  - \int^{0}_{-\gamma \rho} \frac{(-\rho +\gamma\rho^2 + \gamma\rho - 1) + (1 + t)}{(1+t)^2} \,dt \\ &=& (-\rho +\gamma\rho^2 + \gamma\rho - 1)\left( 1 - \frac{1}{1 - \gamma\rho} \right) + \log(1-\gamma\rho) \\ &=&  (1+\rho)\gamma\rho + \log(1-\gamma\rho) =: f(\gamma).
\end{eqnarray*}
To find the optimal step length we need to maximize the function $f$ over $\gamma$. The first order optimality condition gives
\begin{equation*}
    \frac{\partial f}{\partial \gamma}= \frac{\rho^2(- \gamma\rho -  \gamma + 1)}{1 - \rho\gamma} = 0 \quad \Rightarrow \quad 
    \gamma^{\star} = \frac{1}{1 + \rho}.
\end{equation*}
Since $f''(\gamma) = -\frac{\rho^2}{(1 - \gamma\rho)^2} < 0$, we get that $\gamma^{\star}$ is a maximum.

Thus we can solve this problem analytically. The same result was proposed in \cite{nesterov1994interior}, and in \cite{gao2019quasi} it is shown that this step-size maximizes a lower bound on the decrease of the function value. Here we have proved that this step length is actually optimal for this criterion. 

\bibliographystyle{abbrv}
\bibliography{main}

\begin{thebibliography}{1}

\bibitem{burdakov1980some}
O.~P. Burdakov.
\newblock Some globally convergent modifications of {N}ewton's method for
  solving systems of nonlinear equations.
\newblock {\em Doklady Akademii Nauk}, 254(3):521--523, 1980.

\bibitem{de2020worst}
E.~De~Klerk, F.~Glineur, and A.~B. Taylor.
\newblock Worst-case convergence analysis of inexact gradient and {N}ewton
  methods through semidefinite programming performance estimation.
\newblock {\em SIAM Journal on Optimization}, 30(3):2053--2082, 2020.

\bibitem{gao2019quasi}
W.~Gao and D.~Goldfarb.
\newblock Quasi-{N}ewton methods: superlinear convergence without line searches
  for self-concordant functions.
\newblock {\em Optimization Methods and Software}, 34(1):194--217, 2019.

\bibitem{Hildebrand21Newton}
R.~Hildebrand.
\newblock Optimal step length for the {N}ewton method: Case of self-concordant
  functions.
\newblock {\em Math. Methods. Oper. Res.}, 94:253--279, 2021.

\bibitem{nesterov2018lectures}
Y.~Nesterov.
\newblock {\em Lectures on {C}onvex {O}ptimization}, volume 137 of {\em
  Springer Optimization and its Applications}.
\newblock Springer, 2018.

\bibitem{nesterov1994interior}
Y.~Nesterov and A.~Nemirovskii.
\newblock {\em Interior-point polynomial algorithms in convex programming},
  volume~13.
\newblock SIAM, 1994.

\bibitem{pontr}
L.~Pontryagin, V.~Boltyanskii, R.~Gamkrelidze, and E.~Mischchenko.
\newblock {\em The {M}athematical {T}heory of {O}ptimal {P}rocesses.}
\newblock Wiley, New York, London, 1962.

\bibitem{ralph1994global}
D.~Ralph.
\newblock Global convergence of damped {N}ewton's method for nonsmooth
  equations via the path search.
\newblock {\em Mathematics of Operations Research}, 19(2):352--389, 1994.

\end{thebibliography}

\end{document}